\documentclass{amsart}

\usepackage{amssymb}
\usepackage{color}
\usepackage{amsxtra} 
\usepackage{mathrsfs} 
\usepackage{yfonts} 

\newtheorem{theorem}{Theorem}

\newtheorem*{theorem*}{Main Theorem} 
\newtheorem*{theorem**}{Theorem}

\theoremstyle{definition}

\theoremstyle{remark}
\newtheorem{remark}[theorem]{Remark} 
\newtheorem*{remark*}{Remark}

\numberwithin{equation}{section}
\numberwithin{theorem}{section} 

\usepackage{graphicx}

\numberwithin{equation}{section}

\newtheorem*{theorem***}{} 

\theoremstyle{remark} 

\newcommand\restr[2]{{
  \left.\kern-\nulldelimiterspace 
  #1 
  \vphantom{\big|} 
  \right|_{#2} 
}}
\begin{document}

\title[Axisymmetric NS flow 
on the no-slip flat boundary 
]
{
A local instability mechanism of the \\
Navier-Stokes flow with swirl on \\
the no-slip flat boundary
}

\author{Leandro Lichtenfelz}
\address{
Department of Mathematics, University of Pennsylvania, David Rittenhouse Lab., Philadelphia, PA 19104-6395. } 
\email{llichte2@gmail.com} 

\author{Tsuyoshi Yoneda}
\address{Graduate School of Mathematical Sciences, University of Tokyo, Komaba 3-8-1 Meguro, Tokyo 153-8914, Japan} 
\email{yoneda@ms.u-tokyo.ac.jp}

\subjclass[2000]{}

\date{\today} 



\begin{abstract} 
Using numerical simulations of the axisymmetric Navier-Stokes equations with swirl on a no-slip flat boundary, Hsu-Notsu-Yoneda [J. Fluid Mech. 2016] observed the creation of a high-vorticity region on the boundary near the axis of symmetry.
In this paper, using a differential geometric approach, we prove that such flows indeed have a destabilizing effect, which is formulated in terms of a lower bound on the $L^\infty$-norm of derivatives of the velocity field on the boundary.
\end{abstract} 

\maketitle

\section{Introduction} 
\label{sec:Intro}

The Navier-Stokes equations with no-slip flat boundary are expressed as 
\begin{eqnarray}\label{NS}
\partial_tu+(u\cdot \nabla)u-\nu\Delta u+\nabla p&=&0\quad \text{in}\quad \mathbb{R}^3_+\times[0,T),\\
\nonumber
u_0=u|_{t=0},\quad u|_{\partial\mathbb{R}^3}=0,\quad \nabla\cdot u&=&0\quad\text{in} \quad\mathbb{R}^3_+\times[0,T),
\end{eqnarray}
on $\mathbb{R}^3_+ = \{ x = (x_1, x_2, x_3) \in \mathbb{R}^3 : x_3 \geq 0\}$, where $u=u(t, x)$ ($x\in\mathbb{R}^3_+$) is a vector field representing the velocity of the fluid, and $p=p(t, x)$ is the pressure. 

\par The authors in \cite{HNY} performed numerical computations on the axisymmetric Navier-Stokes flow, with and without swirl, on the no-slip flat boundary. They observed a new phenomenon: the distance from the $z$ axis to the points where the velocity attains maximum magnitude exhibits high oscillations (in space) in the swirl case, but not in the swirl-free case. A similar phenomenon happens with the magnitude of vorticity $\|\omega\|$ and vorticity direction $\xi = \omega/\|\omega\|$. These results suggest the sudden appearance of high-vorticity regions near the axis on the boundary. This is what we mean by ``instability" here.

\par Our main purpose in this paper is to prove that the swirl-dominant axisymmetric Navier-Stokes flow with no-slip flat boundary does indeed induce an instability, as observed numerically by Hsu-Notsu-Yoneda \cite{HNY}. Under certain assumptions on the swirl, we get a lower bound on the $L^{\infty}$ norm of derivatives of the velocity field on the boundary.

\par We begin by recalling the aforementioned numerical results of Hsu-Notsu-Yoneda. For the convenience of the reader, we reproduce the relevant figures from \cite{HNY}.
\begin{figure}\label{Fig1}
	\centering
	\includegraphics[scale=0.65]{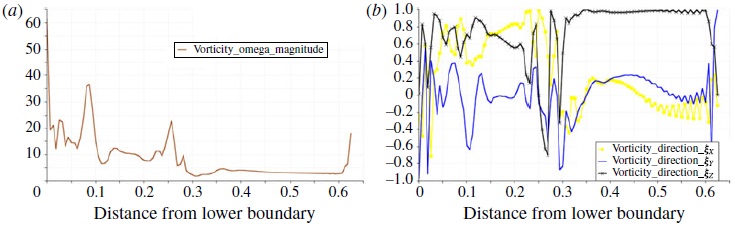}
	\caption{ (From \cite{HNY}) Distance from lower boundary versus: (a) magnitude of vorticity ($\|\omega\|$) and (b) vorticity directions (components of $\xi$) at time $t = 0.4$ on the line parallel to the $z$ axis through the point $(0, 0.05, -0.125)$.}
\end{figure}
In these simulations, the boundary is contained in the plane $z = -0.125$ and therefore meets the axis of symmetry at the point $(0, 0, -0.125)$. At time $t=0.4$, the length of vorticity attains its local maximum at the point $(0, 0.05, -0.125)$ on the boundary, and its vorticity direction is almost parallel to the $x_1$-axis direction (see Figure 2 below: at $x_2=0.05$, $\xi_Y$ is almost zero, while $\xi_X$ is almost $-1.0$. 
Note that the $X$-axis corresponds to the $x_1$-axis and the $Y$-axis corresponds to the $x_2$-axis). 

\begin{figure}\label{Fig2}
	\centering
	\includegraphics[scale=0.65]{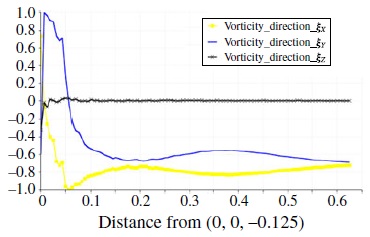}
	\caption{(From \cite{HNY}) Distance from saddle point $(0, 0, -0.125)$ versus vorticity directions (components of $\xi$) at time $t = 0.4$ on the line parallel to the $x_2$-axis through the point $(0, 0.05, -0.125)$.\vspace{-0.1cm}}
\end{figure}
This means that, at the point of maximum vorticity, the swirl component should be dominant compared to
the hyperbolic flow component (we define precisely the swirl and hyperbolic flow components below).

\par We first consider the Navier-Stokes equations in cylindrical coordinates. Write $u=u_re_r+u_\theta e_\theta+u_ze_z$ with $e_r:=(\cos\theta,\sin\theta,0)$, $e_\theta:=(-\sin\theta,\cos\theta,0)$ , $e_z=(0,0,1)$,  $u_r=u_r(t,r,z)$, $u_\theta=u_\theta(t,r,z)$ and  $u_z=u_z(t,r,z)$,
$r=|x|$, $\theta=\arctan (x_2/x_1)$ and $z=x_3$.
Then the axisymmetric Navier-Stokes equations can be expressed as 
\begin{eqnarray}\label{axisymmetricNS}
\nonumber
\partial_t u_r+u_r\partial_ru_r+u_z\partial_zu_r-\frac{u_\theta^2}{r}+\partial_r p&=&\Delta u_r-\frac{u_r}{r^2},\\
\partial_tu_\theta+u_r\partial_ru_\theta+u_z\partial_zu_\theta+\frac{u_ru_\theta}{r}&=&\Delta u_\theta-\frac{u_\theta}{r^2},\\
\nonumber
\partial_tu_z+u_r\partial_ru_z+u_z\partial_zu_z+\partial_zp&=&\Delta u_z,\\
\nonumber
\frac{\partial_r(ru_r)}{r}+\partial_zu_z&=&0
\end{eqnarray}
in $\mathbb{R}^3_+$, where $\Delta=\partial_r^2+(1/r)\partial_r+\partial_z^2$.
The no-slip condition implies 
\begin{equation*}
u_r=u_z=u_\theta=0 \quad\text{on}\quad \partial\mathbb{R}^3_+.
\end{equation*} 
Since the vector fields $\partial_r$ and $\partial_{\theta}$ are tangent to the boundary, we have $\partial_r u = 0$ and $\partial_{\theta} u = 0$ on $\partial\mathbb{R}^3_+$. The divergence free condition then implies $\partial_z u_z = 0$, so that
\begin{equation}\label{eq_vh}
\partial_z u = (\partial_z u_r ) e_r + (\partial_z u_{\theta}) e_{\theta} \quad\text{on}\quad \partial\mathbb{R}^3_+.
\end{equation}

\noindent As we mentioned before, we are interested in the case where the swirl component $\partial_z u_\theta$ is large compared to the radial component $\partial_z u_r$ (this component is corresponding to the hyperbolic flow). Therefore, we assume that $\partial_z u_\theta \neq 0$ except for the origin and introduce the rate of swirl $S$ as follows:
\begin{equation*}
S(t,r)=\frac{\partial_zu_\theta(t,r,0)}{\sqrt{(\partial_zu_r(t,r,0))^2+(\partial_zu_\theta(t,r,0))^2}}.
\end{equation*}

\par Let $v_h$ be the restriction of $\partial_z u$ to the boundary $\partial\mathbb{R}^3_+$. Equation $\eqref{eq_vh}$ shows that $v_h$ is tangent to $\partial\mathbb{R}^3_+$ and since $v_h$ does not vanish it makes sense to consider, for each fixed time $t \in [0, T)$, the integral curve
\begin{gather}\begin{split}
\varphi_t : (-\delta, +\delta) &\rightarrow \partial\mathbb{R}^3_+ \\
\bar \theta &\mapsto \varphi_t(\bar \theta)
\end{split}\end{gather}
of the normalized vector field $v_h(t,x) / |v_h(t,x)|$ starting at 
a reference point $y$ on $\mathbb{R}^2\setminus\{0\}$.
Thus, supressing $t$ we have
$$\partial_{\bar \theta} \varphi(\bar \theta) = \frac{v_h(\varphi(\bar \theta))}{| v_h(\varphi(\bar \theta))|},\quad
\varphi(0)=y\in\mathbb{R}^2\setminus\{0\}
$$
which we can rewrite using the rate of swirl as
$$
\frac{v_h(x)}{| v_h(x)|}= 
S(r) e_\theta(\theta)+\sqrt{1-(S(r))^2} e_r(\theta),
$$
where $\theta=\arctan (x_2/x_1)$ and $r=|x|$.

\begin{remark}
The ``swirl-dominant flow" condition can be expressed as 
\begin{equation}\label{swirl}
	|S(r)|\approx 1\quad\text{and}\quad |\partial_rS(r)|\ll 1.
\end{equation}
in order to have the curvature 
$\kappa:=|\partial_{\bar\theta}^2\varphi|\Big|_{\bar\theta=0}\approx 1/|y|$.
We assume \eqref{swirl} in the main theorem.
\end{remark}

\par To show a local instability of swirl-dominant flows, we define a coordinate system along $\varphi$. For each fixed time $t \in [0, T)$ and reference point $y \in \mathbb{R}^2\backslash\{ 0\}$, the map
\begin{gather}\label{eq_T}
\begin{split}
(z, \bar r, \bar\theta) \mapsto \varphi(\bar \theta) + \frac{\varphi''(\bar \theta)}{|\varphi''(\bar \theta)|}\bar r + ze_z,
\end{split}
\end{gather}
defines a coordinate system near $(y, 0) \in \mathbb{R}^3$, since our assumption on $S$ guarantees that the curvature of $\varphi(\bar \theta)$ does not vanish. We can then express $u$ in terms of a frame relative to this coordinate system as
\begin{equation}
u = u_zdz + u_{\bar r}d\bar r + u_{\bar \theta}fd\bar \theta,
\end{equation}
where $f(\bar r, \bar \theta) = 1 - \bar r\kappa(\bar \theta)$ and $\kappa(\bar \theta) = |\varphi''(\bar \theta)|$ is the curvature of $\varphi(\bar\theta)$.

\par Now, let us formulate ``high-vorticity region near the axis of symmetry" in a mathematical way. It is natural to assume that the high-vorticity region includes a local maximum of $|v_h|$, so let us parametrize the local maxima of $\xi$ as follows: for each $t \in [0, T)$, choose a $\xi(t)$ such that
\begin{equation}\label{eq_xi}
\max_x |v_h(t,x)|=|v_h(t,\xi(t))|.
\end{equation}
Note that because of axisymmetry, equation $\eqref{eq_xi}$ does not uniquely define $\xi(t)$, and in this situation one cannot in general expect it to be everywhere smooth. However, for our purposes it is enough to assume that $\xi(t)$ can be chosen in a $C^1$ way, except for finitely many jump discontinuities.

\par We now give the main theorem which describes the local behavior along $\xi(t)$.
\begin{theorem}
Let $1/M := \sup_{0<t<T} |\xi(t)|$. Assume that $\xi(t)$ is piecewise differentiable and
\begin{equation*}
	\sup_{t\in[0,T), 0<r<M}|S(t,r)|\approx 1\quad\text{and}\quad \sup_{t\in[0,T), 0<r<M}|\partial_rS(t,r)|\ll 1.
\end{equation*}
Define $G$ and $F$ as follows:
\begin{eqnarray*}
F(t,y)&=& 
|\partial_z^3 u_{\bar\theta}(t,z,\bar r,\bar\theta)|+|\partial_z\partial_{\bar r} \partial_{\bar\theta} u_{\bar r}(t,z,\bar r,\bar \theta)|\\
& &
+|\partial_z\partial_{\bar\theta}^2 u_{\bar\theta}(t,z,\bar r,\bar \theta)|+|\partial_z\partial_{\bar r}^2 u_{\bar\theta}(t,z,\bar r,\bar \theta)|\bigg|_{z,\bar r,\bar \theta=0}\\
G(t,y) &=&
\partial_{\bar r}u_{\bar\theta}(t,z,\bar r,\bar\theta)\bigg|_{z,\bar r,\bar \theta=0}
\end{eqnarray*}
Let us take $N>0$ arbitrarily large and fix it. If $|G(0,\xi(0))|>TN$, then at least one of the following two cases must happen:
\begin{equation*}
\sup_{t\in[0,T)}|F(t,\xi(t))|>N\quad\text{or}\quad \sup_{t\in[0,T)}|G(t,\xi(t))|>e^{M^2T}(|G(0,\xi(0))|-TN).
\end{equation*}
\end{theorem}
The above estimate tells us that the swirl dominant flow in the high-vorticity region
on the boundary has a destabilizing effect.

The key idea of the proof comes from Chan-Czubak-Yoneda \cite[Section 2.5]{CCY} and goes back to Ma-Wang  \cite[(3.7)]{MW}.
They investigated 2D separation phenomena using geometric methods, and derived a ``local pressure estimate" on a normal coordinate system in $\bar \theta$, $\bar r$ and $\bar z$ variables. In these coordinates, the tangential and normal derivatives of the scalar function $p$ commute, i.e., 
$\partial_{\bar z}\partial_{\bar \theta}p(\bar \theta, \bar r,\bar z)-\partial_{\bar \theta}\partial_{\bar z}p(\bar \theta,\bar r,\bar z)=0$ and $\partial_{\bar z}\partial_{\bar r}p(\bar \theta, \bar r,\bar z)-\partial_{\bar r}\partial_{\bar z}p(\bar \theta,\bar r,\bar z)=0$ 
 (Lie bracket).
This fundamental observation is the key to extract a local mechanism for pressure estimates.

\section{Proof of the main theorem.}

We first compute the full Navier-Stokes equations in differential form, using the coordinate system introduced in $\eqref{eq_T}$. In this coordinate system, the Euclidean metric (henceforth denoted by $g$) takes the form
\begin{equation*}
g_{11} = 1,~g_{22} = 1,~g_{33} = f^2(\bar r, \bar \theta),
\end{equation*}
where $f(\bar r, \bar \theta) = 1 - \bar r\kappa(\bar \theta)$ and $\kappa(\bar \theta) = |\varphi''(\bar \theta)|$ is the curvature of $\varphi(\bar\theta)$. The indexes $1,2,3$ are corresponding to $z$, $\bar r$, $\bar \theta$ respectively. The non-zero Christoffel symbols are
\begin{gather}
\begin{split}
\Gamma_{33}^2 = -f\partial_{\bar r}f,~\Gamma_{23}^3 = \frac{\partial_{\bar r}f}{f},~
\Gamma_{33}^3 = \frac{\partial_{\bar \theta}f}{f}.
\end{split}
\end{gather}
The volume form is $dV = f\,dz\wedge d\bar r \wedge d\bar \theta$. The hodge dual operator $\star$ acts on $2$-forms as
\begin{gather}
\begin{split}
\star dz \wedge d\bar r = f\,d\bar \theta,~ \star d\bar r \wedge d\bar \theta = \frac{1}{f}\,dz,~\star d\bar\theta \wedge dz = \frac{1}{f}\,d\bar r.
\end{split}
\end{gather}

\begin{remark}
In what follows, we will not distinguish between the vector field $u$ and its corresponding $1$-form, defined by $u^{\flat} = g(u, \cdot)$. Moreover, we write $g(u^{\flat}, e^j)$ as $g(u, e_j)$ for simplicity.
\end{remark}
As before, we expand $u$ in the frame
\begin{align}
u = u_zdz + u_{\bar r}d\bar r + u_{\bar \theta}fd\bar \theta,
\end{align}
where $u_{\bar\theta}=u_{\bar\theta}(t,z,\bar r,\bar \theta)$, $u_{\bar r}=u_{\bar r}(t,z,\bar r,\bar \theta)$, $u_{z}=u_{z}(t,z,\bar r,\bar \theta)$.
\vskip 0.1cm
\noindent \textbf{Divergence}. The divergence of $u$ is given by
\begin{gather}\label{eq_div}
\begin{split}
\mathrm{div}\,u &= \star d\!\star\!u = \star d\Big(u_z f d\bar r \wedge d\bar \theta + u_{\bar r} f d\bar\theta \wedge dz + u_{\bar\theta} dz \wedge d\bar r\Big) \\
&= \star\Big( \partial_z(u_z f) + \partial_{\bar r}(u_{\bar r} f) + \partial_{\bar\theta} u_{\bar\theta} \Big) dz \wedge d\bar r \wedge d\bar\theta \\
&= \frac{1}{f}\Big( \partial_z(u_z f) + \partial_{\bar r}(u_{\bar r} f) + \partial_{\bar\theta} u_{\bar\theta} \Big) = 0.
\end{split}
\end{gather}
From this we get that
\begin{equation}\label{eq:div}
\partial_z u_z f + \partial_{\bar r}(u_{\bar r} f) + \partial_{\bar\theta} u_{\bar\theta} = 0,
\end{equation}
and in particular, on the boundary $\partial \mathbb{R}^3_+$, 
\begin{equation}
\restr{\partial_z u_z}{\partial K}  = -\frac{1}{f}\Big( \restr{\partial_{\bar r}(u_{\bar r} f)}{\partial K} + \restr{\partial_{\bar \theta} u_{\bar\theta}}{\partial K} \Big) = 0,
\end{equation}
since $u$ vanishes on $\partial \mathbb{R}^3_+$ along with its $r$ and $\bar\theta$ derivatives.
\newline\newline \noindent{} \textbf{Laplacian}. We have (since ${\star d \star}u = \mathrm{div}\,u = 0$)
\begin{gather*}
\begin{split}
\Delta u &= {\star d \!\star\! d} u = \\
&\left[\partial_{\bar r}\Big((\partial_z u_{\bar r} - \partial_{\bar r}u_z)f\Big) +
\partial_{\bar\theta}\Big((\partial_z u_{\bar \theta} f - \partial_{\bar\theta} u_z)\frac{1}{f}\Big)
\right]\frac{1}{f} dz ~+ \\
&\left[ \Big(\partial_{\bar\theta}\partial_{\bar r}(u_{\bar\theta} f) - \partial_{\bar\theta}^2 u_{\bar r}\Big)\frac{1}{f} - \frac{\partial_{\bar\theta} f}{f^2}\Big(\partial_{\bar r}(u_{\bar\theta} f) - \partial_{\bar\theta} u_{\bar r}\Big) + \partial_z\Big(\partial_{\bar r}u_z - \partial_z u_{\bar r}\Big)  \right]\frac{1}{f} d\bar r ~+ \\
&\left[ \partial_z\Big(\partial_{\bar\theta} u_z - \partial_z u_{\bar\theta} f\Big)\frac{1}{f} + \Big(\partial_{\bar r}\partial_{\bar\theta} u_{\bar r} - \partial_{\bar r}^2(u_{\bar\theta} f)\Big)\frac{1}{f} + \Big(\partial_{\bar r}(u_{\bar\theta} f) - \partial_{\bar\theta} u_{\bar r}\Big)\frac{\partial_{\bar r} f}{f^2} \right] f d\bar\theta.
\end{split}
\end{gather*}

\noindent We now compute the following quantities in the frame $e^1 = dz$, $e^2 = d\bar{r}$, $e^3 = fd\bar{\theta}$, keeping in mind that $u$ and its $\bar r$ and $\bar\theta$ derivatives vanish on the boundary due to the no slip condition. Note also that for $\bar r = 0$, we have $f = 1$ and $\partial_{\bar \theta} f = 0$.
\begin{gather*}
\begin{split}
\partial_z g(\Delta u, e^2)\big\vert_{\bar r = 0} &= \partial_z g(\Delta u, d\bar r) \\
&= \partial_{\bar r}\partial_{\bar\theta}(\partial_z u_{\bar \theta} f) - \partial_{\bar\theta}^2\partial_z u_{\bar r} + \partial_z^2(\partial_{\bar r} u_z - \partial_z u_{\bar r}). \\ \\
\partial_z g(\Delta u, e^3)\big\vert_{\bar r = 0} &= \partial_z g(\Delta u, f\,d\bar\theta) \\
&= \partial_z^2\Big( \partial_{\bar\theta} u_z - \partial_z u_{\bar\theta} f \Big) + \partial_z \Big(\partial_{\bar r}\partial_{\bar\theta} u_{\bar r} - \partial_{\bar r}^2(u_{\bar\theta} f)\Big) \\
&+ \partial_{\bar r} f\, \partial_z \Big(\partial_{\bar r}(u_{\bar \theta} f) - \partial_{\bar\theta} u_{\bar r}\Big). \\ \\
\partial_{\bar r} g(\Delta u, e^1)\big\vert_{\bar r = 0} &= \partial_{\bar r} g(\Delta u, dz) \\
&= \partial_{\bar r}^2\Big((\partial_z u_{\bar r}) f \Big) + \partial_{\bar r}\partial_{\bar\theta}\partial_z u_{\bar\theta} -\left[ \partial_{\bar r}\Big((\partial_z u_{\bar r})f\Big) +
\partial_{\bar\theta}\partial_z u_{\bar\theta} \right]\partial_{\bar r} f.
\end{split}
\end{gather*}

\begin{gather*}
\begin{split}
\partial_{\bar\theta} g(\Delta u, e^1)\big\vert_{\bar r = 0} &= \partial_{\bar\theta} g(\Delta u, dz) \\
&= \partial_{\bar\theta}\partial_{\bar r}\Big((\partial_z u_{\bar r}) f \Big) + \partial_{\bar\theta}^2\partial_z u_{\bar\theta}.
\end{split}
\end{gather*}
\newline \noindent{} \textbf{Covariant Derivative $\nabla_u u$}. From the Christoffel symbols, we see that the only non-zero covariant derivatives between coordinate vector fields are
\begin{gather}
\begin{split}
\nabla_{\partial_{\bar r}}\partial_{\bar\theta} = \nabla_{\partial_{\bar\theta}}\partial_{\bar r} = \frac{\partial_{\bar r} f}{f}\partial_{\bar\theta},~\nabla_{\partial_{\bar \theta}}\partial_{\bar\theta} = f\,\partial_{\bar r} f\, \partial_{\bar r} + \frac{\partial_{\bar\theta} f}{f}\partial_{\bar\theta}.
\end{split}
\end{gather}
It follows that
\begin{gather}
\begin{split}
\nabla_{\partial_z} u &= (\partial_z u_z)\partial_z + (\partial_z u_{\bar r})\partial_{\bar r} + \left(\frac{\partial_z u_{\bar\theta}}{f}\right)\partial_{\bar\theta}. \\
\nabla_{\partial_{\bar r}} u &= (\partial_{\bar r} u_z)\partial_z + (\partial_{\bar r} u_{\bar r}) \partial_{\bar r} + \left(\frac{\partial_{\bar r} u_{\bar \theta}}{f}\right)\partial_{\bar\theta}. \\
\nabla_{\partial_{\bar\theta}} u &= (\partial_{\bar\theta} u_z)\partial_z + \left(\partial_{\bar\theta} u_{\bar r} + u_{\bar r}\frac{\partial_{\bar r} f}{f} + \partial_{\bar r} f u_{\bar\theta} \right)\partial_{\bar r} + \left( \frac{\partial_{\bar\theta} u_{\bar\theta}}{f} \right)\partial_{\bar\theta}.
\end{split}
\end{gather}
The tangential components are
\begin{gather}
\begin{split}
g(\nabla_u u, e_2) &= g(\nabla_u u, \partial_{\bar r}) \\
&= u_r \partial_z u_{\bar r} + u_{\bar r} \partial_{\bar r} u_{\bar r} + u_{\bar\theta}\left(\partial_{\bar\theta} u_{\bar r} + u_{\bar r}\frac{\partial_{\bar r} f}{f} + \partial_{\bar r} f u_{\bar\theta} \right). \\ \\
g(\nabla_u u, e_3) &= \frac{1}{f}g(\nabla_u u, \partial_{\bar\theta}) \\
&= \frac{1}{f^2}\left( u_z \partial_z u_{\bar\theta} + u_{\bar r} \partial_{\bar r}u_{\bar\theta} + u_{\bar\theta}\partial_{\bar\theta} u_{\bar\theta} \right)\,g(\partial_{\bar\theta}, \partial_{\bar\theta}) \\
&= u_z \partial_z u_{\bar\theta} + u_{\bar r} \partial_{\bar r}u_{\bar\theta} + u_{\bar\theta}\partial_{\bar\theta} u_{\bar\theta}.
\end{split}
\end{gather}
We now take a derivative in the direction orthogonal to the boundary ($z$ direction), and evaluate everything at $\bar r=0$.
\begin{gather}
\begin{split}
\partial_z g(\nabla_u u, e_2)\big\vert_{\bar r = 0} &= \partial_z u_z\partial_z u_{\bar r}. \\
\partial_z g(\nabla_u u, e_3)\big\vert_{\bar r = 0} &= \partial_z u_z\partial_z u_{\bar\theta}.
\end{split}
\end{gather}
\newline \noindent{} \textbf{Pressure ($\nabla p$)}. In order to eliminate the pressure and have an equation in terms of $u$ alone, we want to use the fact that the no-slip condition gives $dp = \Delta u$ on the boundary. Thus, we write
\begin{gather}
\begin{split}
\partial_z g(\nabla p, e_3)\big\vert_{\bar r = 0} &= \frac{\partial_z \partial_{\bar\theta} p}{f(0)} = \partial_{\bar\theta} g(\nabla p, e_1) = \partial_{\bar\theta} g(dp, e^1) = \partial_{\bar\theta} g(\Delta u, e^1) \\
&= \partial_{\bar\theta} \Big[ \partial_{\bar r}\big( (\partial_z u_{\bar r}) f \big) + \partial_{\bar\theta}\partial_z u_{\bar\theta}  \Big]. \\ \\
\partial_z g(\nabla p, e_2)\big\vert_{\bar r = 0} &= \partial_z \partial_{\bar r} p = \partial_{\bar r} g(\nabla p, e_1) = \partial_{\bar r} g(\Delta u, e_1) \\
&= \partial_{\bar r}^2\Big((\partial_z u_{\bar r}) f \Big) + \partial_{\bar r}\partial_{\bar\theta}\partial_z u_{\bar\theta} -\left[ \partial_{\bar r}\Big((\partial_z u_{\bar r})f\Big) +
\partial_{\bar\theta}\partial_z u_{\bar\theta} \right]\partial_{\bar r} f.
\end{split}
\end{gather}
\newline \noindent{} \textbf{ODE for $\partial_z u_{\bar\theta}$}. Using the above formulas we can derive an ODE for $\partial_z u_{\bar\theta}$ at $(\bar r, \bar \theta, z) = (0, 0, 0)$. Let $t_0 > 0$ be some positive time. Then, we have
\begin{gather}\label{eq:integral}
\begin{split}
\partial_z u_{\bar\theta}(t_0, 0) &= \partial_z g(u(t_0), e_3) = \partial_z g(u_0, e_3) + \int\limits_0^{t_0} \frac{d}{dt} \partial_z g(u, e_3)\, dt \\
&= \int\limits_0^{t_0} \partial_z g(\partial_t u, e_3) + g(u, \partial_t e_3)\,  dt \\
&= \int\limits_0^{t_0} \partial_z g(\Delta u, e_3) -\partial_z g(\nabla_{\! u} u, e_3) -\partial_z g(\nabla p, e_3) + g(u, \partial_t e_3)\, dt.
\end{split}
\end{gather}
We expand the integrand above and write it in terms of the quantities:
\begin{gather}
\begin{split}
\alpha_1(t,y) = \partial_z u_{\bar\theta}(t, z,\bar r,\bar\theta)\bigg|_{z,\bar r,\bar \theta=0}, \hspace{0.5cm} \beta_1(t,y) &= \partial_z\partial_{\bar r} u_{\bar\theta}(t, z,\bar r,\bar \theta)\bigg|_{z,\bar r,\bar \theta=0}, \\
\alpha_2(t,y) = \partial_z u_{\bar r}(t, z,\bar r,\bar\theta)\bigg|_{z,\bar r,\bar \theta=0}, \hspace{0.5cm} \beta_2(t,y) &= \partial_z\partial_{\bar\theta} u_{\bar r}(t, z,\bar r,\bar\theta)\bigg|_{z,\bar r,\bar \theta=0}, \\
\alpha_3(t,y) = \partial_z^3 u_{\bar\theta}(t, z,\bar r,\bar\theta)\bigg|_{z,\bar r,\bar \theta=0}, \hspace{0.5cm} \beta_3(t,y) &= \partial_z\partial_{\bar r} \partial_{\bar\theta} u_{\bar r}(t, z,\bar r,\bar\theta)\bigg|_{z,\bar r,\bar \theta=0}, \\
\beta_4(t,y) &= \partial_z\partial_{\bar\theta}^2 u_{\bar\theta}(t,z,\bar r,\bar\theta)\bigg|_{z,\bar r,\bar \theta=0}, \\
\beta_5(t,y) &= \partial_z\partial_{\bar r}^2 u_{\bar\theta}(t,z,\bar r,\bar \theta)\bigg|_{z,\bar r,\bar \theta=0}.
\end{split}
\end{gather}

\noindent Recall that $\partial_z u_z = 0$ by the divergence free condition and that $f(0) = 1$, $\partial_{\bar r} f(0) = -\kappa(0)$, where $\kappa(\bar\theta)$ is the curvature of $\varphi(\bar\theta)$. We abbreviate $\kappa = \kappa(0)$, $\partial_{\bar\theta}\kappa = \partial_{\bar\theta} \kappa(0)$. The Laplacian term becomes
\begin{gather}
\begin{split}
\partial_z g(\Delta u, e_3) &= \partial_z^2\Big( \partial_{\bar\theta} u_z - \partial_z u_{\bar\theta} f \Big) + \partial_z \Big(\partial_{\bar r}\partial_{\bar\theta}u_{\bar r} - \partial_{\bar r}^2(u_{\bar\theta} f)\Big) + \partial_{\bar r} f \partial_z \big(\partial_{\bar r} (u_{\bar\theta} f) - \partial_{\bar\theta} u_{\bar r}\big) \\ \\
&= \partial_z\partial_{\bar\theta}(\partial_z u_z) - \alpha_3 + \beta_3 - \partial_{\bar r}^2(\partial_z u_{\bar\theta} f) -\kappa \partial_{\bar r}(\partial_z u_{\bar\theta} f) + \kappa\beta_2 \\ \\
&= \partial_z\partial_{\bar\theta}(\partial_z u_z) - \alpha_3 + \beta_3 - \beta_5 + 2\kappa\beta_1 -\kappa\beta_1 + \kappa^2\alpha_1 + \kappa\beta_2 \\ \\
&= \partial_z\partial_{\bar\theta}(\partial_z u_z) - \alpha_3 + \beta_3 - \beta_5 + \kappa\beta_1 + \kappa^2\alpha_1 + \kappa\beta_2.
\end{split}
\end{gather}
At this point, we use the divergence equation $(\ref{eq:div})$ to rewrite $\partial_z u_z$ in terms of $\partial_{\bar r} u_{\bar r}$ and $\partial_{\bar\theta} u_{\bar\theta}$ instead, which leads to
\begin{gather}
\begin{split}
\partial_z\partial_{\bar\theta}(\partial_z u_z) = -\beta_3 + \partial_{\bar\theta}\kappa\, \alpha_2 + \kappa\beta_2 - \beta_4.
\end{split}
\end{gather}
Putting everything together, we get
\begin{gather}
\begin{split}
\partial_z g(\Delta u, e_3) &= -\beta_3 + \partial_{\bar\theta}\kappa \,\alpha_2 + \kappa\beta_2 - \beta_4  - \alpha_3 \\
&+ \beta_3 - \beta_5 + \kappa\beta_1 + \kappa^2\alpha_1 + \kappa\beta_2 \\
&= \partial_{\bar\theta}\kappa\, \alpha_2 + 2\kappa\beta_2 - \beta_4  - \alpha_3 -\beta_5 + \kappa\beta_1 + \kappa^2\alpha_1.
\end{split}
\end{gather}
Similarly, one can check that
\begin{gather}
\begin{split}
\partial_z g(\nabla p, e_3) = -\kappa'\alpha_2 -\kappa \beta_2 + \beta_3 + \beta_4,
\end{split}
\end{gather}
and $\partial_zg(\nabla_u u, e_3)(0) = 0$. Thus,
\begin{gather}
\begin{split}
&\int\limits_0^{t_0} \partial_z g(\Delta u, e_3) -\partial_z g(\nabla_{\! u} u, e_3) -\partial_z g(\nabla p, e_3)\,  dt \\
= &\int\limits_0^{t_0} 2\partial_{\bar\theta}\kappa\, \alpha_2 - \kappa\beta_2 -2\beta_4 -\beta_3 -\alpha_3 -\beta_5 + \kappa\beta_1 + \kappa^2\alpha_1  \,dt.
\end{split}
\end{gather}
It follows that $\alpha_1(t,y)$
 satisfies the ODE
\begin{equation}\label{eq:ODE}
\partial_t \alpha_1 = \kappa^2\alpha_1 + (2\partial_s\kappa + \rho  )\alpha_2 -\alpha_3 +\kappa\beta_1 -\kappa\beta_2 -\beta_3 -2\beta_4 -\beta_5.
\end{equation}

\noindent Note that $\alpha_2=0$, and  along $\xi$,   
\begin{eqnarray*}
& &\partial_t\left(\alpha_1(t,\xi(t))\right)=
\partial_t\alpha_1(t,y)+\nabla_y\alpha_1(t,y)\bigg|_{y=\xi(t)}\cdot \partial_t\xi(t)=
\partial_t\alpha_1(t,y),\\
& &\beta_1(t,\xi(t))=\beta_2(t,\xi(t))=0.
\end{eqnarray*}
Thus
\begin{equation*}
\partial_t (\alpha_1\circ\xi) = \kappa^2\alpha_1\circ\xi -\left(\alpha_3+\beta_3 +2\beta_4+\beta_5\right)\circ\xi.
\end{equation*}
Let $F=-\alpha_3-\beta_3-2\beta_4-\beta_5$,
and assume $|F|<N$.
In our case we have 
\begin{equation}
\begin{split}
(\alpha_1\circ\xi)(t) &= \int_0^t e^{^{
{
\int_{\tau}^t \kappa^2(y)dy}}} (F\circ\xi)(\tau)d\tau
+ (\alpha_1\circ\xi)(0)e^{^{{ \int_{0}^t \kappa^2(y)dy}}} \\
&\geq
(\alpha_1\circ\xi)(0)e^{^{{ \int_{0}^t \kappa^2(y)dy}}}  -N\int_0^t e^{^{{ \int_{\tau}^t \kappa^2(y)dy}}}d\tau\\
&\geq
e^{^{{ \int_{0}^t \kappa^2(y)dy}}}\left((\alpha_1\circ\xi)(0)  -N\int_0^t e^{^{{ -\int_{0}^{\tau} \kappa^2(y)dy}}}d\tau\right)\\
&\geq
e^{M^2T}\left((\alpha_1\circ\xi)(0)  -NT\right).
\end{split}
\end{equation}
This is the desired estimate.

\vspace{0.5cm}
\noindent
{\bf Acknowledgments.} The authors would like to thank Doctor Pen-Yuan Hsu 
for helpful comments.
TY  was partially supported by
Grant-in-Aid for Young Scientists A (17H04825)
and Scientific Research B (18H01135), 
Japan Society for the Promotion of Science (JSPS).


\bibliographystyle{amsplain}

\end{document}